\documentclass[12pt]{article}

\usepackage{amsmath,amssymb}

\newtheorem{theorem}{Theorem}
\newtheorem{lemma}[theorem]{Lemma}

\begin{document}

\title
 {Movable algebraic singularities of second-order ordinary differential equations}

\author{G. Filipuk$^1$ and R. G. Halburd$^{2,1}$}

\maketitle

\begin{tabular}{ll}
$\!\!\!\! ^1\!$ Department of Mathematical Sciences\ \ \ \ \ \ \ \ \ \ \ \ \ \ \ 	& $\!\!\!\! ^2\!$ Department of Mathematics\\
Loughborough University				& University College London\\
Loughborough					 	& Gower Street \\
Leicestershire LE11 3TU				& London WC1E 6BT\\
UK								& UK\\
\  & \\ 
{\tt G.Filipuk@lboro.ac.uk}			& {\tt R.Halburd@ucl.ac.uk}
 \end{tabular}



\begin{abstract}
Any nonlinear equation of the form $y''=\sum_{n=0}^N a_n(z)y^n$ has a (generally branched) solution with leading order behaviour proportional to $(z-z_0)^{-2/(N-1)}$ about a point $z_0$, where the coefficients $a_n$ are analytic at $z_0$ and  $a_N(z_0)\ne 0$.  We consider the subclass of equations for which each possible leading order term of this form corresponds to a one-parameter family of solutions represented near $z_0$ by a Laurent series in fractional powers of $z-z_0$.  For this class of equations we show that the only movable singularities that can be reached by analytic continuation along finite-length curves are of the algebraic type just described.  This work generalizes previous results of S. Shimomura.  The only other possible kind of movable singularity that might occur is an accumulation point of algebraic singularities that can be reached by analytic continuation along infinitely long paths ending at a finite point in the complex plane.  This behaviour cannot occur for constant coefficient equations in the class considered.  However, an example of R. A. Smith shows that such singularities do occur in solutions of a simple autonomous  second-order differential equation outside the class we consider here.
\end{abstract}

\section{Introduction}
The most fundamental result in the theory of differential equations in the complex domain is Cauchy's existence and uniqueness theorem.
\begin{theorem} (Cauchy)
Let $f_1,\ldots,f_n$ be analytic functions in a neighbourhood of the point
$(z_0,\eta_1,\ldots,\eta_n)$ in $\mathbb{C}^{n+1}$.  There is a unique $n$-tuple of functions 
$y_1,\ldots,y_n$ analytic in a neighbourhood $\Omega$ of $z_0\in\mathbb{C}$ such that
for all $z\in\Omega$,
$$
y'_i=f_i(z;y_1,\ldots,y_n)\quad\mbox{   and   }\quad  y_i(z_0)=\eta_i,\quad\mbox{for all}\  i=1,\ldots,n.
$$
\end{theorem}

Cauchy's theorem is a \textit{local} result.  It guarantees the existence of a local solution of a differential equation provided that the equation is well behaved (analytic) at the initial condition.  However, it says nothing about the nature of this solution after analytic continuation.  For example, Cauchy's theorem guarantees that
the initial value problem 
$$
y'=\frac1{2(z+1)}\left(y-y^3\right),\qquad y(0)=c,
$$
has a unique solution in a neighbourhood of $z=0$.  This solution is 
$$
y(z)=c[(1+z)/(1+c^2z)]^{1/2}.
$$
For $c\ne 0$ this solution is analytic on $\{z\,:\,|z|<\min\{1,|c|^{-2}\}$ but has  branch points at $z=-c^{-2}$ and $z=-1$.  The location of the singularity at $z=-c^{-2}$ varies with the initial condition.  Such singularities are called \textit{movable}.
The location of  the singularity at $z=-1$ does not depend on the initial condition.  Such singularities are called \textit{fixed}.  Heuristically speaking, fixed singularities occur at points where the equation itself is singular in some sense, whereas movable singulaities develop spontaneously.  
All singularities of autonomous equations are movable.

A singularity of a solution of the linear differential equation
$$
y^{(n)}+a_{n-1}(z)y^{(n-1)}+\cdots+a_1(z)y'+a_0(z)y=f(z)
$$
can only occur at a singularity of one of the coefficient functions $a_0,\ldots,a_{n-1}$ or $f$.  It follows that linear equations do not have movable singularities.   In \cite{painleve:88}, Painlev\'e showed that the only movable singularities of solutions of first-order ODEs of the form
\begin{equation}
\label{Palgebraic}
F(z;y,y')=0,
\end{equation}
where $F$ is a polynomial in $y$ and $y'$  with coefficients that are analytic in some common domain,
are poles or algebraic branch points.  The situation for higher-order equations is much more complicated and very poorly understood.  For example, the general solution of the equation
$$
(yy''-y'^2)^2+4yy'^3=0
$$
is $y(z)=c\exp\{(z-z_0)^{-1}\}$, which has a movable essential singularity at $z_0$.  
Painlev\'e considered the equation
\begin{equation}
\label{painleve}
y''=\frac{2y-1}{y^2+1}y'^2,
\end{equation}
which has the general solution 
\begin{equation}
\label{painsoln}
y(z)=\tan(\log(c_1z-c_2)),
\end{equation}
where $c_1$ and $c_2$ are constants.
If $c_1\ne 0$, this solution has movable poles accumulating at the movable branched singularity at $z=c_2/c_1$.  Painlev\'e's example is often used as a warning about the overzealous use of so-called Painlev\'e tests.  Using formal series methods, it is very easy to show that equation \eqref{painleve} has a two-parameter family of Laurent series solutions about poles.  Although every nontrivial solution contains infinitely many movable poles, it also contains a branched accumulation point of poles, which is not detected by the standard analysis.

Of particular relevance to our investigation is the example of Smith
\cite{smith:53},
\begin{equation}
\label{smith}
y''+4y^3y'+y=0.
\end{equation}
He showed that if a solution of \eqref{smith} can be continued analytically along a curve $\gamma$ of finite length up to but not including the point $z_0$, then $y$ has an algebraic branch point at $z_0$.  Moreover, he showed that there are solutions such that analytic continuation along a curve of infinite length ends at a point in the finite $z$-plane that is an accumulation point of these algebraic branch points.   This kind of accumulation of singularities is quite different from that in Painlev\'e's example.  We will discuss this example in more detail in section \ref{sect-accumulate}. 

For third order equations other phenomena arise.  For example, the general solution of the Chazy equation
$$
y'''=2yy''-3y'^2
$$
has a movable natural barrier \cite{chazy:09,chazy:10}.  

In this paper we will primarily be concerned with movable singularities of solutions of the equation
\begin{equation}
\label{original}
y''(z)=\sum_{n=0}^Na_n(z)y(z)^n,\qquad n\in\mathbb{N},\quad n\ge 2.
\end{equation}
A singularity at $z_0$ of a solution of equation \eqref{original} is called movable if $a_0,\ldots,a_N$ are analytic at $z_0$ and $a_N(z_0)\ne 0$.
Note that the case $0\le N\le 1$ is linear, in which case $y$ will have no movable singularities.  For the rest of this paper we assume that $N\ge 2$.
The main result of this paper is the following theorem, which is a characterization of the movable singularities of equation \eqref{original}.  It is a generalization of results of Shimomura described below.

\begin{theorem}
\label{main}
For $N\ge 2$, suppose that  there is a domain $\Omega\subset\mathbb{C}$ such that $a_0,\ldots,a_N$ are analytic and that $a_N(z_0)\ne 0$ on $\Omega$.  Suppose further that for each $z_0\in\Omega$ and for each $c_0$ such that
\begin{equation}
\label{c0def}
c_0^{N-1}=\frac{2}{a_N(z_0)}\frac{N+1}{(N-1)^2},
\end{equation}
equation \eqref{original}
admits a formal series solution of the form 
\begin{equation}
\label{origseries}
y(z)=\sum_{j=0}^\infty c_j(z-z_0)^{\frac{j-2}{N-1}}.
\end{equation}
Then
\begin{enumerate}
\item
For each $c_0$ satisfying \eqref{c0def} and for each $\beta\in\mathbb{C}$, there is a unique formal series solution of the form \eqref{origseries} such that $c_{2(N+1)}=\beta$.  
\item
Given $c_0$ and $c_{2(N+1)}$ as above, the series \eqref{origseries} converges in a neighbourhood of $z_0$.
\item
Now let $y$ be a solution of equation \eqref{original} that can be continued analytically along a curve 
$\gamma$ up to but not including the endpoint $z_0$, where the coefficients $a_j$ are analytic on $\gamma\cup\{z_0\}$ and $a_N$ is nowhere zero on $\gamma\cup\{z_0\}$.  If $\gamma$ is of finite length, then $y$ has a convergent series expansion about $z_0$ of the form \eqref{origseries}.  
\item
If $y$ cannot be represented by a series expansion about $z_0$ of the form \eqref{origseries} then $\gamma$ is of infinite length and $z_0$ is an accumulation point of such algebraic singularities.
\end{enumerate}
\end{theorem}

By ``accumulation point'' in part {4} of the theorem we mean that given any $\epsilon>0$ there exists a straight line segment $l$ in the disk of radius $\epsilon$ centred at $z_0$ with endpoints $z_1\in\gamma$ and $z_2$ such that analytic continuation of $y$ along $\gamma$ up to $z_1$ and then along $l$ ends in an algebraic singularity at $z_2$.

We will see in section 2.3 that the assumptions of theorem \ref{main} correspond to one or two differential relations relating the coefficients of equation \eqref{original} depending on whether $N$ is even or odd respectively.   These relations are resonance conditions.  That is, they represent the vanishing of an obstruction that occurs when the recurrence relation determining the coefficients in the series 
\eqref{origseries} breaks down.  We will see in section 2.1 that equation \eqref{original} can be normalized so that $a_N=2(N+1)/(N-1)^2$ and $a_{N-1}=0$.   With this normalization the resonance condition in the even $N$ case is $a''_{N-2}=0$.  In the odd $N$ case the two resonance conditions are equivalent to $a''_{N-2}=0$ and another condition, which is listed below for the first few odd values of  $N$:
\begin{gather*}
N=3:\quad  a_0'(z)=0,\\
N=5:\quad  \left[4a_1(z)-a^2_3(z)\right]' =0,\\
N=7:\quad \left[10a_2(z)-9a_4(z)a_5(z)\right]'=0.
\end{gather*}

The fact that the resonance conditions are satisfied if each $a_n$ is a constant is apparent from the following ``explicit'' integration. If $y$ is  non-constant then multiplying  equation \eqref{original} by $2y'$ and integrating gives
\begin{equation}
\label{firstintegral}
y'^2=\kappa+\sum_{n=0}^N \frac{2a_n}{n+1}y^{n+1},
\end{equation}
for some constant $\kappa$.  Equation \eqref{firstintegral}, which can be solved in terms of hyperelliptic functions or their degenerations, is of the form \eqref{Palgebraic}.  Therefore, we know that the only singularities of its solutions are algebraic, so the accumulation of singularities described by part 4 of the theorem does not occur here.  It remains an open question whether such an accumulation of singularities can occur for equations of the form \eqref{original}.  This phenomenon is known to occur in the example of Smith 
\eqref{smith}, which lies outside of the class considered here.

The role of movable singularities is particularly important in the theory of integrable systems in the context of the Painlev\'e property.
An ODE is said to possess the Painlev\'e property if all movable singularities of all solutions are 
poles\footnote{Some authors say that an ODE is said to possess the Painlev\'e property if all solutions are single-valued around all movable singularities.}.  If solutions of an equation possess no fixed singularities in the finite plane, e.g., if the equation is autonomous, then the Painlev\'e property is equivalent to the property that all solutions are meromorphic.  The connection between the integrability of an ODE and the singularity structure of its solutions appears to have been first exploited by Kowalevskaya, who used it to find a new integrable case of the equations of motion of a spinning top \cite{kowalevskaya:89a,kowalevskaya:89b}.

In \cite{fuchs:84}, L.~Fuchs studied equations of the form 
$$
y'=F(z;y),
$$
where $F$ is rational in $y$ and with coefficients that are analytic in some domain.
He showed that the only nonlinear equation in this class with the Painlev\'e property is the Riccati equation, which corresponds to the case in which $F$ is quadratic in $y$.  The general solution of the Riccati equation
\begin{equation}
\label{riccati}
y'=a(z)y^2+b(z)y+c(z),
\end{equation}
where $a\not\equiv 0$,
is given by
\begin{equation}
\label{riccatitrans}
y(z)=-\frac{1}{a(z)}\frac{w'(z)}{w(z)},
\end{equation}
where $w$ is the general solution of the linear equation
$$
a(z)\frac{d^2w}{dz^2}-\left[a'(z)+a(z)b(z)\right]\frac{dw}{dz}
+c(z)a^2(z)w=0.
$$
Since the singularities of $w$ are fixed, we see from the transformation
(\ref{riccatitrans}) that the only movable singularities of $y$ are poles, corresponding to zeros of $w$. 
Here we see our first example of simple movable singularities (poles) detecting integrable (in this case, linearizable) equations.

Painlev\'e, Gambier and Fuchs classified all equations with the Painlev\'e property of the form 
$$
y''=F(z;y,y'),
$$
where $F$ is rational in $y$ and $y'$ with coefficients that are analytic in some domain.  They found that each such equation could be transformed to one of fifty canonical forms.  Six of these canonical equations are now called the \textit{Painlev\'e equations} $P_I$--$P_{\,V\!I}$.  All of the other canonical equations could either be solved in terms of classically known functions, by quadrature, by solving linear ODEs or in terms of solutions of one of $P_I$--$P_{\,V\!I}$.   The Painlev\'e equations were subsequently shown to be compatibility conditions for certain isomonodromy (linear) problems.  
They are extremely important equations in the theory of integrable systems.  The solutions of the Painlev\'e equations are often referred to as nonlinear special functions.  For each Painlev\'e equation, all solutions are meromorphic on the universal cover of $\mathbb{C}\setminus S$, where $S$ is a finite (possibly empty) subset of $\mathbb{C}$.

Interest in the Painlev\'e property and movable singularities was rekindled in the late twentieth century through connections with the spin-spin correlation function of the two-dimensional Ising model (Wu, McCoy, Tracy and Barouch \cite{wmtb:76,mccoytw:77}),
holonomic quantum fields (Sato, Miwa and Jimbo \cite{satomj:78}) and the theory of soliton equations. 
In their studies of the asymptotic behaviour of solutions of soliton equations
Ablowitz and
Segur~\cite{ablowitzs:77}, and Ablowitz, Ramani, and
Segur~\cite{ablowitzrs:78,ablowitzrs:80} discovered that ODE reductions of soliton equations often led to equations of Painlev\'e type, possibly after a
transformation of variables.  

Relatively simple tests of strong necessary conditions for the Painlev\'e property can be used to detect integrability
\cite{ablowitzrs:78,ablowitzrs:80}.  These tests involve showing that formal Laurent series solutions exist and that certain resonance conditions are satisfied.    However, it is a much more difficult task to show that a given equation actually possesses the Painlev\'e property, without having a more-or-less explicit representation of the solution.

The main idea of many of the standard proofs of the Painlev\'e property is to characterize the singular Laurent series expansions in terms of locally analytic variables satisfying a regular initial value problem.  The main tool in this respect is the following result, which can be found in a number of sources \cite{gromakls:02,hille:76} and follows from the fact that Cauchy's existence and uniqueness theorem can be strengthened to give a lower bound on the radius of convergence of the solution.

\begin{lemma} (Painlev\'e)
\label{regularsequence}
Let $f_1,\ldots,f_m$ be analytic functions in a neighbourhood of the point
$(\alpha,\eta_1,\ldots,\eta_m)$ in $\mathbb{C}^{m+1}$. 
Let $\gamma$ be a curve with end point $\alpha$ and suppose that $y_i$ is analytic on 
$\gamma\setminus\{\alpha\}$ for $i=1,\ldots,m$
and satisfies
$$
y'_i=f_i(z;y_1,\ldots,y_m).
$$
Let $(z_n)$ be a sequence of points such that $z_n\in\gamma$, $z_n\to\alpha$ and
$y_i(z_n)\to\eta_i$ as $n\to\infty$ , for all $i=1,\ldots,n$.  Then each $y_i$ is analytic at $\alpha$.
\end{lemma}

Painlev\'e himself provided a proof that the first Painlev\'e equation
\begin{equation}
\label{p1}
y''=6y^2+z
\end{equation}
possesses the Painlev\'e property.  This proof, which appears in a number of forms in the literature
(e.g., Ince \cite{ince:56}, Golubev \cite{golubev:50}), had a number of gaps in it that have been filled by several authors (Hukuhara \cite{hukuhara:60,okamotot:01}, Hinkkanen and Laine \cite{hinkkanenl:99}, Shimomura \cite{shimomura:03}).  All of these proofs proceed in the following manner.
\begin{itemize}
\item[\textit{i.}]
Show that if a solution $y$ of equation \eqref{p1} has a pole at $z=z_0$ then it is a double pole and the coefficient $\beta$ of $(z-z_0)^4$ is arbitrary.  Fixing a value for $\beta$ uniquely determines all the other coefficients.
\item[\textit{ii.}]
Consider a finite length curve on which $y$ is analytic except at the endpoint $z_0$ where it is singular.  Show that $y$ is unbounded on $\gamma$.  Show that if $y$ is not bounded away from zero on $\gamma$ then $\gamma$ can be deformed to a new finite length curve $\tilde\gamma$ ending at $z_0$ so that $y$ is bounded away from zero on $\tilde\gamma$.
\item[\textit{iii.}]
Using the form of the series expansion in \textit{i}, introduce a new variable $u:=y^{-1/2}$, for some choice of branch and another variable $v$ such that $v(z_0)$ encodes the value of the resonance parameter $\beta$.  The functions $u$ and $v$ satisfy a regular initial value problem with $u(z_0)=0$ and $v(z_0)=\kappa$ for $\kappa\in \mathbb{C}$.
\item[\textit{iv.}]
A function $W$ of $z$, $y$ and $y'$ is shown to be bounded on $\gamma$.  This function is in some sense an approximate first integral of the equation.
\item[\textit{v.}]
If $A:=\liminf_{\gamma\ni z\to\infty}|y(z)|$ is finite and positive then the boundedness of $W$ leads to $y_1=y$ and $y_2=y'$ satisfying the conditions of lemma \ref{regularsequence}, so $y$ is analytic at $z_0$.
\item[\textit{vi.}]
If $A=\infty$ (i.e. if $\lim_{\gamma\ni z\to\infty}|y(z)|=\infty$), then the boundedness of $W$ (and the correct choice of branch in the definition of $u$) shows that $v$ must be bounded on $\gamma$, so applying  lemma \ref{regularsequence} to the initial value problem for $u$ and $v$ mentioned in \textit{iii} shows that $y$ must have a double pole at $z_0$.
\end{itemize}

Many of the other methods for proving that the Painlev\'e equations possess the Painlev\'e property rely on the underlying isomonodromy problems and Riemann-Hilbert techniques (Miwa \cite{miwa:81}, Malgrange \cite{malgrange:82}, see also Fokas, Its, Kapaev and Novokshenov  \cite{fokasikn:06}).  These methods explicitly exploit structure that is closely associated with the integrability of the Painlev\'e equations and they are therefore unlikely to generalize directly to the nonintegrable equations that we consider.  Other proofs of the Painlev\'e property for the Painlev\'e equations are Steinmetz \cite{steinmetz:00}, who uses differential inequalities, Erugin \cite{erugin:76}, and Joshi and Kruskal \cite{joshik:94}.

Shimomura \cite{shimomura:05,shimomura:07}  studied the quasi-Painlev\'e property, which is also known as the weak Painlev\'e property (see \cite{ramanidg:82,ramanigb:89}).  An ODE is said to possess the quasi-Painlev\'e property if all movable singularities of all solutions are at most algebraic branch points.  He proved that for any $k\in\mathbb{N}$, analytic continuation of any solution $y$ of 
\begin{equation}
\label{shimomura}
y''=\frac{2(2k+1)}{(2k-1)^2}y^{2k}+z
\end{equation}
along a finite length curve ends at a point where either $y$ is analytic, has a pole or has an algebraic branch point.   The method of proof is essentially the same as outlined above for equation \eqref{p1} except that the initial value problem for the new variables $u$ and $v$ described in \textit{iii} is now for $z(u)$ and $v(u)$, rather than $u(z)$ and $v(z)$.  In this way it is shown that $z$ is an analytic function of $u$ but its inverse is in general analytic in a fractional power of $z-z_0$.  

The main result of the present paper is a generalization of Shimomura's result.  Due to the complexity and generality of the series expansions we would otherwise have to consider, we show directly that the existence of algebraic formal series solutions is equivalent to the existence of a bounded function $W$ (modulo the same curve modification arguments outlined in the description of the proof of the Painlev\'e property in \textit{ii}).  The choice of suitable functions $u$ and $v$ comes directly from the expression for $W$.  The equations for $W$, $u$ and $v$ are implicit in the sense that an algorithm is presented for calculating their coefficients in terms of the functions $a_j$ and that these algorithms are well defined provided that the resonance conditions are satisfied.  This is especially useful in the odd $N$ case in which one of the resonance conditions is determined recursively.  It is this use of representing certain associated equations using recurrence relations that allows us to get a results for such a large class of equations.  The rest of the proof follows the general pattern above.

The possibility that algebraic branch points could accumulate along infinite length curves in a bounded region of the complex plane appeared in Smith's study of equations of the form 
$$
y''+f(y)y'+g(y)=P(z),
$$
where $f$ and $g$ are polynomials and the degree of $f$ is greater than the degree of $g$ and $P$ is analytic in some domain \cite{smith:53}.

In section 2 we will transform equation \eqref{original} to a canonical form that will simplify future calculations.  We will then show how to construct singular series expansions of solutions.  The different nature of the leading orders and the resonance conditions in the even $N$ and odd $N$ cases will be described.  In section 3 we show that the existence of these singular series expansions implies the boundedness of a certain function $W$.  This section  then follows closely the outline above.  In section 4 we consider analytic continuation along infinite length curves in a bounded region. We end with a discussion in section 5.

\section{Formal series expansions}

In this section we use standard methods to determine when formal Laurent series expansions of solutions of equation (\ref{original}) in fractional powers of $z-z_0$ exist.  The method described in this section is the most direct way of verifying whether a given equation of the form \eqref{original} satisfies the assumptions of theorem \ref{main}.  The fundamental difference between the even $N$ and odd $N$ cases will be immediately apparent.  Although this section only deals with formal series expansions, it will be shown at the end of section 3 that they converge.

\subsection{Canonical form of equation \eqref{original}}
Suppose that $a_0,\ldots,a_N$ are analytic in a neighbourhood of $z_0$ and that $a_N(z_0)\ne 0$. 
We begin by transforming equation \eqref{original} to a canonical form in order to simplify our analysis.  To this end,  let
$$
f(z)=\left(
	\frac{2(N+1)}{(N-1)^2a_N(z)}
	\right)^{1/(N+3)},
\quad
g(z)=-\frac{a_{N-1}(z)}{Na_N(z)};
$$
$$
y(z)=f(z)\tilde y(\tilde z)+g(z),
\quad
\tilde z=\int_{z_0}^z f^{-2}(\tau)\,{\rm d}\tau.
$$
This transformation is analytic and invertible in a neighbourhood of $z=z_0$.  Under this transformation, $z=z_0$ is mapped to $\tilde z=0$ and equation \eqref{original} is mapped to an equation of the same form with $z$ and $y(z)$ replaced by $\tilde z$ and $\tilde y(\tilde z)$ and with $a_j(z)$ replaced by
$\tilde a_j(\tilde z)$, $j=0,\ldots, N$,  where all the coefficients $\tilde a_j$ are analytic in  a neighbourhood $\tilde\Omega$ of $0$,
$\tilde a_{N}(\tilde z)=2(N+1)/(N-1)^2$ and $\tilde a_{N-1}(\tilde z)=0$.
Throughout sections 2 and 3 we will replace equation \eqref{original} by this canonical form.  By a slight abuse of notation, we will not use tildes on the new variables but replace them by the orginals.
Hence, without loss of generality, we restrict our attention to the equation
\begin{equation}
\label{canonical}
y''(z)=\sum_{n=0}^{N-2}a_n(z)y(z)+\frac{2(N+1)}{(N-1)^2}y^N(z),
\end{equation}
in an open set $\Omega$ in which the coefficients $a_n$ are analytic.

\subsection{Leading order analysis}
Choose $\alpha\in\Omega$.  We will look for formal solutions of equation \eqref{canonical} of the form
\begin{equation}
\label{formalpq}
y(z)=\sum_{n=0}^\infty c_n\zeta^{np-q},
\end{equation}
where $p$ and $q$ are positive, $c_0\ne 0$ and $\zeta=z-\alpha$.  Substituting the expansion \eqref{formalpq} into equation \eqref{canonical} and keeping only the leading-order terms (i.e., the terms that are largest for small $\zeta$) on both sides, we have
$$
q(q+1)c_0\zeta^{-q-2}+\dots=\frac{2(N+1)}{(N-1)^2}c_0^N\zeta^{-qN}+\cdots,\quad \zeta\to 0.
$$
Equating powers of $\zeta$ in the leading order terms gives
$q=2/(N-1)$.  Equating the coefficients of
$\zeta^{-2N/(N-1)}$ then gives $c_0^{N-1}=1$.  So the leading order behaviour of $y$ is
$$
y\sim c_0\zeta^{-2/(N-1)},\qquad c_0^{N-1}=1.
$$
In terms of the original variables in equation \eqref{original}, this condition on $c_0$ corresponds to equation \eqref{c0def}.

It is at this point that we first see that the parity of $N$ plays an important role.  If $N$ is even then the mapping $\zeta\mapsto\zeta^{-2/(N-1)}$  has $N-1$ branches about $\zeta=0$.  The different choices of $c_0$ correspond to different choices of branch.   The choice of $c_0$ can be effectively absorbed into the choice of branch, so we can take $c_0=1$, i.e., $y(z)\sim \zeta^{-2/(N-1)}$.  If $N=2K+1$ is odd then the mapping
$\zeta\mapsto\zeta^{-2/(N-1)}=\zeta^{-1/K}$ has $K=(N-1)/2$ branches, whereas there are $N-1=2K$ choices of $c_0$.   
For some $\epsilon>0$, let $x_0$ be one of the possible values of the function $f(z)=z^{-1/K}$ at $z=\epsilon$.  Analytic continuation of $f$ in a clockwise direction $n$ times around the circle $|z|=\epsilon$ gives the value $\omega^nx_0$, where $\omega=\exp(2\pi {\rm i}/K)$ is the primitive $K{\mbox{th}}$ root of unity.
The choices for $c_0$ naturally fall into two classes: those for which $c_0^K=1$ and those for which $c_0^K=-1$.  For a fixed determination of $\zeta^{-1/K}$, it can be seen that analytic continuation around $\zeta=0$ cannot change the effective value of $c_0$ from one class to the other because analytic continuation effectively only multiplies $c_0$ by a power of $\omega$.  Hence, absorbing the choice of $c_0$ as much as possible into the choice of branch, we see that the odd $N$ case has two essentially different possible leading order behaviours.

\subsection{Resonance conditions\label{sect-resonance}}
Having computed the possible leading order behaviours, we now consider how they can be extended to series solutions.  A natural guess would be to extend them as Laurent series in $\zeta^{2/(N-1)}$. In this case the left side of equation \eqref{canonical} would also be a Laurent series in $\zeta^{2/(N-1)}$.  However, $\zeta$ can be expanded as a Laurent series in $\zeta^{2/(N-1)}$ about $0$ if and only if $N$ is odd (in which case $\zeta=(\zeta^{2/(N-1)})^{(N-1)/2}$).  Therefore, the right side of equation \eqref{canonical} would also be a Laurent series in $\zeta^{2/(N-1)}$ if $N$ is odd, but not if $N$ is even and some of the coefficient functions $a_0,\ldots,a_{N-2}$ are not constant.  However, if $N$ is even and we extend the leading order term to a Laurent series in
$\zeta^{1/(N-1)}$, then both sides of equation \eqref{canonical} will also be Laurent series in this variable.

\vskip 5mm

\noindent{\bf Proof of theorem \ref{main} part 1}

We look for formal series solutions of equation \eqref{canonical} of the form
\begin{equation}
\label{generalalg}
y(z)=\sum_{n=0}^\infty c_n\zeta^{\frac{n-2}{N-1}},
\end{equation}
which corresponds to \eqref{formalpq} with $p=1/(N-1)$ and $q=2/(N-1)$.
For even $N$ we have $c_0=1$.  For odd $N$ we have $c^{(N-1)/2}_0=\pm 1$ and $c_{2k+1}=0$ for all $k=0,1,2,\ldots$  Substituting \eqref{generalalg} into equation \eqref{canonical} and equating coefficients of like powers of $\zeta$ gives a recurrence relation for the coefficients $c_n$.  For any positive integer $r$, the lowest power of $\zeta$ such that $c_r$ appears as a coefficient on either side of equation \eqref{canonical} is $q_r:=(r-2N)/(N-1)$.  The coefficient of $\zeta^{q_r}$ on the left side of \eqref{canonical} is
$$
\frac{(r-2)(r-N-1)}{(N-1)^2}c_r,
$$
while the  corresponding coefficient on the right side has the form
$$
\frac{2N(N+1)}{(N-1)^2}c_r+P_r(c_0,\ldots,c_{r-1}),
$$
where $P_r$ is a polynomial in its arguments.  Equating these coefficients gives a recurrence relation of the form
\begin{equation}
\label{usualrec}
(r+N-1)(r-2N-2)c_r=(N-1)^2P_r(c_0,\ldots,c_{r-1}),
\end{equation}
for all $r=1,2,\ldots$

Equation \eqref{usualrec} is a recurrence relation for the coefficients $c_r$ with a \textit{resonance} at
$r=2(N+1)$.  In other words, we can use \eqref{usualrec} to determine $c_r$ inductively in terms of $c_0,\ldots,c_{r-1}$ except when $r=2(N+1)$.  A necessary and sufficient condition for the existence of a formal series expansion of the form \eqref{generalalg} is that $P_{2(N+1)}(c_0,\ldots,c_{r-1})=0$.  If this condition is satisfied then $c_{2(N+1)}$ is arbitrary.  
\hfill$\Box$

\vskip 5mm

In the even $N$ case there is effectively just one leading order behaviour, so we get just one resonance condition.  We will see later that this condition is $a_{N-2}''(\alpha)=0$.   Since $\alpha\in \Omega$ is arbitrary, this becomes the differential condition
$a_{N-2}''\equiv 0$.  That is $a_{N-2}(z)=Az+B$ for some constants $A$ and $B$.  In the odd $N$ case, there are two leading order behaviours leading to two  resonance conditions.  These two conditions taken together lead to the condition $a_{N-2}''=0$ and another condition, which we express in a more convenient way (in equations \eqref{recurrence} and \eqref{secondconstraint}, which are still not in closed form)  in the next section.


\section{Movable algebraic branch points}

We begin this section by establishing the properties of certain functions that will be useful in the proof of part 3 of theorem \ref{main}.   Part 2 of theorem \ref{main} will be a consequence of the existence of a particular regular initial value problem that arises in the proof of part 3 of theorem \ref{main}.

\subsection{A bounded function}
Let $y$ be a solution of equation \eqref{canonical} that can be analytically continued along a finite length curve $\gamma$ up to but not including the endpoint $z_0$.  The purpose of this section is to identify a function $W(z;y,y')$ that is bounded along $\gamma\setminus\{z_0\}$ whenever $y$ is bounded away from 0 on $\gamma\setminus\{z_0\}$.  If the coefficients $a_n$ are all constants then a suitable function would be the first integral 
$$
y'^2-2\sum_{k=1}^{N+1}\frac{a_{k-1}}{k}y^k.
$$
More generally, we consider the ansatz
\begin{equation}
\label{Wdef}
W(z):=y'(z)^2+\left(\sum_{k=1}^{N-1}\frac{b_k(z)}{y^k(z)}\right)y'(z)-2\sum_{k=1}^{N+1}\frac{a_{k-1}(z)}{k}y^k(z),
\end{equation}
where the functions $b_1,\ldots, b_{N-1}$ are analytic at $z_0$ but yet to be determined.

We will now obtain an integral representation for $W$.
Differentiating equation (\ref{Wdef}) with respect to $z$ and substituting the second derivative of $y(z)$ from equation  (\ref{canonical})  and $y'^2$ from (\ref{Wdef}), we get a linear equation of the form
\begin{equation}\label{WDiff}
W'+P(z,1/y)W=Q(z,1/y)y'+R(z,1/y)+S(z,y),
\end{equation}
where 
\begin{gather}\nonumber
P(z,1/y)=\sum_{k=1}^{N-1}kb_k(z)/y^{k+1},\\\nonumber
Q(z,1/y)=\sum_{k=1}^{N-1}b_k'(z)/y^k+\left(\sum_{k=1}^{N-1}kb_k(z)/y^{k+1}
\right)\left(\sum_{k=1}^{N-1}b_k(z)/y^k\right),\\ \label{RSterms} 
\mbox{and}\qquad R(z,1/y)+S(z,y)=\\   \nonumber
\sum_{k=1}^{N-1}b_k(z)/y^k\sum_{j=0}^Na_j(z)y^j
-2
\sum_{k=1}^{N+1}a_{k-1}'(z)y^k/k
-2
\left(
\sum_{k=1}^{N-1}k b_k(z)/y^{k+1}
\right)
\left(
\sum_{k=1}^{N+1}a_{k-1}y^k/k
\right).
\end{gather}
The functions $P,\;Q,\;R,\;S$ are polynomials in their second arguments. We separate non-positive powers of $y(z)$ into $R(z,1/y) $ and positive into $S(z,y).$ Equating the coefficients of the positive powers of $y(z)$ in \eqref{RSterms} gives
\begin{equation}
\label{Ssimple}
S(z,y)=\sum_{n=1}^{N-1}\bigg(\bigg(\sum_{m=1}^{N-n}\frac{n-m+1}{n+m+1}b_m(z)a_{m+n}(z)\bigg)-\frac{2}{n}a_{n-1}'(z)\bigg)y^n.
\end{equation}

Now consider the case in which $N$ is even and use the recurrence relation
\begin{equation}
\label{recurrence}
\frac{N+1-2n}{(N-1)^2}b_n(z)=\frac{1}{N-n}a_{N-n-1}'(z)-\frac{1}{2}\sum_{m=1}^{n-1}\frac{N-n-m+1}{N-n+m+1}b_{m}(z)a_{N+m-n}(z),
\end{equation}
to define
uniquely $b_1,\ldots,b_{N-1}$, so that $S(z,y)\equiv0$.  Note that the coefficient of $b_n(z)$ on the left side of equation \eqref{recurrence} does not vanish for $n=1,\ldots,N-1$ since $N$ is even.

For any point $z$ on $\gamma$, we denote by $\gamma_z$ the part of $\gamma$ up to $z$
 (in particular, $\gamma=\gamma_{z_0}$).  Using the integrating factor
\begin{equation}
\label{Edef}
E(z):=\exp\left(
\int_{\gamma_z}P(\zeta,1/y(\zeta))\,{\rm d}\zeta
\right),
\end{equation}
the solution of equation \eqref{WDiff} can be written as
\begin{equation}
\label{integrated}
W(z)=\frac{1}{E(z)}
\left(
\kappa+
\int_{\gamma_z}
\left\{
Q(\tau,1/y(\tau))y'+R(\tau,1/y(\tau))
\right\}
E(\tau)\,{\rm d}\tau
\right),
\end{equation}
for some constant $\kappa$.  

We now impose the main assumption of theorem \ref{main}, namely that there is a formal series solution $y$ of equation \eqref{canonical} that is a Laurent series in $(z-z_0)^{1/(N-1)}$ with leading order behaviour given by 
$y\sim (z-z_0)^{-2/(N-1)}$.  For such a solution, $E(z)$, $Q(z,1/y(z))$ and $R(z,1/y(z))$ are power series in $(z-z_0)^{1/(N-1)}$.  Furthermore, $E(z)=1+O\left((z-z_0)^{\frac{N+3}{N-1}}\right)$ and 
$Q(z,1/y(z))y'=-\frac{2b_1'(z_0)}{N-1}(z-z_0)^{-1}+\cdots$.  Hence the integral in equation \eqref{integrated} shows that the series expansion for $W$ about $z_0$ contains a logarithm if $b_1'(z_0)\ne 0$.  However, from the definition \eqref{Wdef}, $W$ has a Laurent series expansion in $(z-z_0)^{1/(N-1)}$.  So $b_1'(z_0)= 0$.  Furthermore, this analysis can be repeated for any point $z_0\in \Omega$, therefore $b_1'\equiv 0$.  On substituting $n=1$ in the recurrence relation \eqref{recurrence}, we see that $b_1=a_{N-2}'$.  So we have proved the necessity of the condition $a''_{N-2}\equiv 0$.

Next we consider the case in which $N=2K+1$ is odd.  We define the first $K=(N-1)/2$ functions
$b_1,\ldots,b_K$ using the recurrence relation \eqref{recurrence}.  From equation 
\eqref{Ssimple}, we see that $S(z,y)$ is a polynomial in $y$ of degree at most $K=(N-1)/2$.  The coefficient of $y^K$ in $S(z,y)$ is
$$
\rho(z)=
\bigg(\sum_{m=1}^{K}\frac{K+1-m}{K+1+m}b_m(z)a_{m+K}(z)\bigg)-\frac{2}{K}a_{K-1}'(z).
$$
By assumption there are two families of solutions at $z_0$ that are Laurent series in $(z-z_0)^{1/K}$ with asymptotic behaviour
$y\sim c_0(z-z_0)^{-1/K}$.  The two families correspond to the cases $c^K_0=1$ and $c^K_0=-1$.  For such solutions,
$P$, $R$ and $E$ are all power series in $(z-z_0)^{1/K}$, where $E$ is again an integrating factor given by
\eqref{Edef}.  We also have the following leading-order behaviours:
$S(z,y(z))\sim\rho c_0^K(z-z_0)^{-1}$, 
$Q(z,1/y(z))y'\sim-\frac{2b_1'(z_0)}{N-1}(z-z_0)^{-1}$, $E(z)\sim 1$.
Integrating equation \eqref{WDiff} we have
\begin{equation}
\label{odd-integrated}
W(z)=\frac{1}{E(z)}
\left(
\kappa+
\int_{\gamma_z}
\left\{
Q(\tau,1/y(\tau))y'+R(\tau,1/y(\tau))+S(\tau,y(\tau))
\right\}
E(\tau)\,{\rm d}\tau
\right).
\end{equation}
The coefficient of $(z-z_0)^{-1}$ in the integrand is now
$$
l(z_0;c_0^K):=\rho(z_0) c_0^K-\frac{2b_1'(z_0)}{N-1}.
$$
Hence to avoid a logarithm in the series expansion of $W$ in both the cases $c^K_0=1$ and $c^K_0=-1$, we must impose the conditions $l(z_0;1)=l(z_0;-1)=0$.  These conditions  are equivalent to $b_1'(z_0)=0$ and $\rho(z_0)=0$ and must hold for all $z_0$ in a nonempty open set, leading to $a''_{N-2}\equiv 0$ and 
\begin{equation}
\label{secondconstraint}
\frac{1}{K}a_{K-1}'(z)-\frac{1}{2}\sum_{m=1}^{K}\frac{K+1-m}{K+1+m}b_m(z)a_{K+m}(z) =0.
\end{equation}
Furthermore, knowing that equation \eqref{secondconstraint} is satisfied, we can now choose $b_{K+1}$ arbitrarily (e.g., $b_{K+1}\equiv 0$) and use the recurrence relation to define the functions $b_k$ for $k=K+2,\ldots, N-1$, which means that $S(z,y)\equiv 0$.

We have shown in both the even $N$ and odd $N$ cases that $W$ defined by \eqref{Wdef}, where
the functions $b_k$ are defined by the recurrence relation \eqref{recurrence} and $b_1'\equiv a_{N-2}''\equiv 0$, is given by
equation \eqref{integrated} for some constant $\kappa$.  Now let $\gamma$ be any  curve in $\Omega$ of finite length $L$ ending at $z_0\in \Omega$ and let $y$ be any function analytic on $\gamma\setminus\{z_0\}$ such that there exists a constant $B\ge 1$ and 
$|y(z)|^{-1}\le B$ for all $z\in\gamma\setminus\{z_0\}$.  Note that we are no longer assuming anything about any possible singularity of $y$ at $z_0$.

Since $\gamma$ and its endpoints are in $\Omega$, there is a constant $A$ and functions
{$p_2,\ldots,p_N;$}\allowbreak{$q_2,\ldots,q_{2N-1};$}\allowbreak{$r_0,\ldots,r_{N-1}$} that are all analytic and bounded by $A$ on $\gamma$
such that
$$
P(z,1/y)=\sum_{k=2}^{N} p_k/y^k,
\quad
Q(z,1/y)=\sum_{k=2}^{2N-1} q_k/y^k
\quad\mbox{and  }
R(z,1/y)=\sum_{k=0}^{N-1} r_k/y^k.
$$
So on $\gamma$,
$$
|P(z,1/y(z))|
\le
\sum_{k=2}^{N} |p_k|/|y|^k
\le 
\sum_{k=2}^{N} AB^k
\le (N-1)AB^{N}
$$
and similarly
$$
|R(z,1/y(z))|
\le NAB^{N-1}.
$$
So from equation \eqref{Edef} and the inequalities $\exp(-|x|)\le |\exp x|\le \exp(|x|)$, we have
\begin{equation}
\label{Eineq}
\exp(-(N-1)AB^{N}L)\le |E(z)|\le\exp((N-1)AB^{N}L).
\end{equation}
Therefore it follows from equation \eqref{integrated} that in order to prove that $W$ is bounded on $\gamma$, it is sufficient to show that
$$
I:=\int_{\gamma_z}
Q(\tau,1/y(\tau))y'
E(\tau)\,{\rm d}\tau
$$
is bounded.  Now
\begin{eqnarray*}
I(z)&=&-\sum_{k=2}^{2N-1}\int_{\gamma_z} \frac{q_k(\tau)E(\tau)}{k-1}\left(\frac{1}{y^{k-1}(\tau)}\right)'\,{\rm d}\tau
\\
&=&
C-\sum_{k=2}^{2N-1} \frac{q_k(z)E(z)}{(k-1)y^{k-1}(z)}
\\
&&
+\sum_{k=2}^{2N-1}\int_{\gamma_z} \frac{q'_k(\tau)+q_k(\tau)P(\tau,1/y(\tau))}{k-1}\frac{E(\tau)}{y^{k-1}(\tau)}\,{\rm d}\tau,
\end{eqnarray*}
for some constant $C$.
The integrand is bounded on $\gamma$, therefore  $I$ is bounded and hence so is $W$.

In summary, we have proved the following.

\begin{lemma}
\label{Wlemma}
For $N\ge 2$, suppose that  there is a domain $\Omega\subset\mathbb{C}$ such that $a_0,\ldots,a_{N-2}$ are analytic on $\Omega$.  Suppose further that for each $z_0\in\Omega$ and for each $c_0$ such that $c_0^{N-1}=1$
equation \eqref{canonical}
admits a formal series solution of the form 
$$
y(z)=\sum_{j=0}^\infty c_j(z-z_0)^{\frac{j-2}{N-1}}.
$$
Then
\begin{equation}
\label{aN2}
a''_{N-2}\equiv 0
\end{equation}
and there are functions $(b_k)_{k=1}^{N-1}$ that are analytic on $\Omega$ and satisfy the recurrence relation \eqref{recurrence}
for $n=1,\ldots,N-1$.
Furthermore, let $\gamma$ be a finite-length curve in $\Omega$ and let $y$ be a solution of equation \eqref{canonical} that is analytic and bounded away from 0 on $\gamma$, then the function 
$W$ defined by equation \eqref{Wdef}
is bounded on $\gamma$.
\end{lemma}

\vskip 5mm

\noindent If $N$ is even then the recurrence relation (\ref{recurrence}) defines the functions $b_k$ uniquely without further constraints.  However, if $N$ is odd then the left side of (\ref{recurrence}) vanishes when $n=(N+1)/2$.  This gives the second constraint \eqref{secondconstraint}.
If this constraint is satisfied, then the function $b_{(N+1)/2}$ is not determined by the recurrence relation
(\ref{recurrence}) and can be chosen to have any value (e.g., zero).  In section \ref{sect-resonance} we saw that the existence of formal Laurent series expansions in fractional powers of $z-z_0$ was equivalent to one or two resonance conditions being satisfied in the even $N$ and odd $N$ cases respectively.  Lemma \ref{Wlemma} shows that equation (\ref{aN2}) (in the even and odd $N$ cases) and equation (\ref{secondconstraint}) (in the odd $N$ case) are necessary consequences of these resonance conditions.  We will see subsequently that these conditions are also sufficient.

\subsection{Proof of theorem \ref{main} part 3}
Let $y$ be a solution of equation (\ref{canonical}) that can be continued analytically along a finite length curve 
$\gamma$ up to but not including the endpoint $z_0$, where the coefficients $a_j$ are analytic on $\gamma\cup\{z_0\}$.  
If $y$ is bounded on $\gamma$ then so is 
$$
y'=C+\int_\gamma \sum_{n=0}^Na_n(\zeta)y(\zeta)^n \,{\rm d}\zeta. 
$$
Hence, from lemma \ref{regularsequence}, $y$ is analytic at $z_0$, which is a contradiction.
So we have shown that $y$ is unbounded on $\gamma$.  In other words, 
$\limsup_{\gamma\ni z\to\infty}|y(z)|=\infty$.
Following the main ideas behind Painlev\'e's work on the Painlev\'e property, we divide the proof into three parts depending on the value of
$$
A:=\liminf_{\gamma\ni z\to\infty}|y(z)|.
$$
\vskip 3mm

\noindent{\bf Case 1:} $0<A<\infty$.\\
Since $y$ is bounded away from $0$ on $\gamma$ near $z_0$, lemma \ref{Wlemma} shows that
$W$ is bounded.  By assumption, there exists a sequence of points $z_n\in\gamma$ such that $z_n\to z_0$ and $y(z_n)\to\eta$ for some $\eta\in\mathbb{C}$.  Viewing equation \eqref{Wdef} as a quadratic equation in $y'$, we see that there is a subsequence of $(z_n)$ on which $y'$ also approaches a finite value.  Hence, applying lemma \ref{regularsequence} to equation \eqref{canonical} with $y_1=y$ and $y_2=y'$, we see that $y$ is regular at $z_0$, which is a contradiction.

\vskip 3mm

\noindent{\bf Case 2:} $A=\infty$.\\
In this case $y\to\infty$ as $z\to 0$ on $\gamma$ and $W$ is bounded on $\gamma$ near $z_0$.  Solving equation \eqref{Wdef} for $y'$ gives
\begin{equation}
\label{firstyprime}
y'=-\frac 12
\left\{
	\left(
		\sum_{k=1}^{N-1}\frac{b_k(z)}{y^k(z)}
	\right)
	-
	\left[
		\left(
		\sum_{k=1}^{N-1}\frac{b_k(z)}{y^k(z)}
		\right)^2
		+8
		\sum_{k=1}^{N+1}\frac{a_{k-1}(z)y^k(z)}{k}
		+4W(z)
	\right]^{1/2}
\right\},
\end{equation}
for some choice of branch of the square root.  

Assume for now that $N$ is even.  We introduce a new variable $u$ such that $y=1/u^2$.  There is a choice of branch if we take this as a ``definition'' of $u$, which we will fix shortly.  Define the function
$F(z;u)$, which is analytic for $u$ near zero and $z$ near $z_0$, by
\begin{eqnarray}
F(z;u)^2
&=&
u^{2(N+1)}\left\{
			\left(
		\sum_{k=1}^{N-1}{b_k(z)}{u^{2k}}
		\right)^2
		+8
		\sum_{k=1}^{N+1}\frac{a_{k-1}(z)}{ku^{2k}}
		\right\}
		\nonumber
		\\
&=&
\frac{16}{(N-1)^2}
+
8\sum_{k=1}^N\frac{a_{k-1}(z)}{k}u^{2(N-k+1)}
+
\sum_{k,l=1}^{N-1} b_k(z)b_l(z)u^{2(N+k+l+1)},
\label{Fsqr}
\end{eqnarray}
where the choice of branch is determined such that $F(z,0)=4/(N-1)$.  Now equation \eqref{firstyprime}
takes the form
\begin{equation}
\label{anotheryprime}
y'=-\frac12\sum_{k=1}^{N-1}b_ku^{2k}+\frac{F(z;u)G(z;u,W)}{2u^{N+1}},
\end{equation}
where
$$
G(z;u,W)=\left\{
		1+\frac{4u^{2(N+1)}W}{F(z;u)^2}
		\right\}^{1/2}=1+\frac{2u^{2(N+1)}W}{F(z;u)^2}+\cdots
$$
and the choice of the ``$+$'' sign in equation \eqref{anotheryprime} determines the branch in the definition of $u$ (essentially $u^{N+1}=y^{-(N+1)/2}$).  Now we introduce one more variable $v$ such that the expansions for $G$ in $W$ and $v$ agree to first order.  That is, we define $v$ such that
\begin{equation}
\label{evenvdef}
G(z;u,W)=\left\{
		1+\frac{4u^{2(N+1)}W}{F(z;u)^2}
		\right\}^{1/2}=1+\frac{2u^{2(N+1)}v}{F(z;u)^2},
\end{equation}
giving
\begin{equation}
\label{evenW}
W=v+\frac{u^{2(N+1)}}{F(z;u)^2}v^2.
\end{equation}
Since $y\to\infty$ as $z$ approaches $z_0$ along $\gamma$, it follows that $u\to 0$ along $\gamma$.
From the boundedness of $W$ and the definition of $v$ \eqref{evenvdef}, it follows that  $v$ is bounded  on $\gamma$.  We will now construct an initial value problem for $z$ and $v$ as functions of $u$ that is regular at $u=0$.  

Equation \eqref{anotheryprime} takes the form
$$
y'=-\frac12\sum_{k=1}^{N-1}b_ku^{2k}+\frac{F(z;u)}{2u^{N+1}}+\frac{u^{N+1}}{F(z;u)}v.
$$
Recalling that $y=u^{-2}$, we have
\begin{equation}
\label{evenueqn}
\frac{dz}{du}=u^{-(N-2)}J(z;u,v),
\end{equation}
where $J$ is analytic in its arguments at $(z_0;0,\kappa)$, for any $\kappa\in\mathbb{C}$, and
$J(z_0;0,\kappa)=1-N$.  

Now we will construct an equation for $v'$.  Differentiating equation \eqref{evenW} with respect to $z$ gives
\begin{equation}
\label{evenWdiff}
W'=\left(1+\frac{2u^{2(N+1)}}{F(z;u)^2}v\right)v'+\left(\frac{u^{2(N+1)}}{F(z;u)^2}\right)'v^2.
\end{equation}
Dividing equation \eqref{Fsqr} by $u^{2(N+1)}$, differentiating with respect to $z$ and using the reciprocal of equation \eqref{evenueqn} to replace $du/dz$ on the right side, shows that
\begin{equation}
\label{evenfracest}
\left(\frac{u^{2(N+1)}}{F(z;u)^2}\right)'=K(z;u,v),
\end{equation}
where $K$ is analytic in its arguments around $(z;u,v)=(z_0;0,\kappa)$, for any $\kappa\in\mathbb{C}$.
Finally, from equation \eqref{WDiff} with $S(z,y)\equiv 0$, we have
\begin{eqnarray}
W'&=&-P(z,u^2)
	\left(
		v+\frac{u^{2(N+1)}}{F(z;u)^2}v^2
	\right)
	-2
	Q(z,u^2)
	u^{-(N+1)}J(z;u,v)
	+R(z,u^2)
	\nonumber
\\
&&
=
u^{-(N-3)}L(z;u,v),
\label{evenWest}
\end{eqnarray}
where $L$ is analytic in its arguments around $(z;u,v)=(z_0;0,\kappa)$, for any $\kappa\in\mathbb{C}$. 
Using equations \eqref{evenfracest} and \eqref{evenWest} in equation \eqref{evenWdiff} shows that
$v'$ is given by an equation of the form
$$
\frac{dv}{dz}=u^{-(N-3)}H(z;u,v),
$$
where $H$ is analytic in its arguments around $(z;u,v)=(z_0;0,\kappa)$, for any $\kappa\in\mathbb{C}$. 
Hence
\begin{equation}
\label{evenveqn}
\frac{dv}{du}=\frac{dv}{dz}\frac{du}{dz}=uJ(z;u,v)H(z;u,v).
\end{equation}

Note that the system \eqref{evenueqn} and \eqref{evenveqn} is regular at $(z;u,v)=(z_0;0,\kappa)$. Since $u\to 0$ and $v$ is bounded on $\gamma$ near $z_0$, it follows from Painlev\'e's lemma (lemma \ref{regularsequence}) that $z$ and $v$ are analytic  functions of $u$ in a neighbourhood of $u=0$.  Recall that $J(z_0;0,\kappa)\allowbreak =1-N$. Hence from equation \eqref{evenueqn}, we see that $z$ has a series expansion in $u$ of the form
$$
z=z_0+\sum_{n=0}^\infty \alpha_nu^{n+N-1},
$$
where $\alpha_0=1$.  Subtracting $z_0$ and taking the $(N-1)$st root gives
$$
(z-z_0)^{1/(N-1)}=u\left(1+\sum_{n=1}^\infty\beta_n u^n\right),
$$
for some choice of branch on the left side.
Inverting this series shows that $u$ is a power series in $(z-1)^{1/(N-1)}$ of the form 
$$
u=\sum_{n=1}^\infty \lambda_n (z-z_0)^{n/(N-1)},
$$
where $\lambda_0=1$.  It follows that $y$ has a series expansion of the required form.

Next we consider the case in which $N=2K+1$ is odd.  Returning to equation \eqref{firstyprime}, we  let $u=1/y$, giving
\begin{equation}
\label{oddanotheryprime}
y'=-\frac12\sum_{k=1}^{N-1}b_ku^{k}-\varepsilon\frac{F(z;u)G(z;u,W)}{2u^{K+1}},
\end{equation}
where $\varepsilon=\pm 1$, $F(z,u)$ is defined near $(z_0,0)$ by
\begin{equation}
F(z;u)^2
=
\frac{4}{K^2}
+
8\sum_{k=1}^{2K}\frac{a_{k-1}(z)}{k}u^{2K-k+2}
+
\sum_{k,l=1}^{2K} b_k(z)b_l(z)u^{2K+k+l+2},
\label{oddFsqr}
\end{equation}
$F(z,0)=2/K$ and
$$
G(z;u,W)=\left\{
		1+\frac{4u^{2(K+1)}W}{F(z;u)^2}
		\right\}^{1/2}=1+\frac{2u^{2(K+1)}W}{F(z;u)^2}+\cdots
$$
Note that the $\varepsilon=\pm 1$ in equation \eqref{oddanotheryprime} cannot be absorbed into a suitable choice of branch without affecting the value of $F$ or $G$ at $u=0$.  Following a similar argument as in the even $N$ case, we define $v$ by
$$
G(z;u,W)=\left\{
		1+\frac{4u^{2(K+1)}W}{F(z;u)^2}
		\right\}^{1/2}=1+\frac{2u^{2(K+1)}v}{F(z;u)^2}.
$$
We then find that $z$ and $v$ satisfy a system of differential equations of the form
\begin{eqnarray}
\frac{dz}{du} &=& u^{K-1}A(z;u,v),\label{oddueqn}\\
\frac{dv}{du} &=& B(z;u,v),\label{oddveqn}
\end{eqnarray}
where $A$ and $B$ are analytic in $(z;u,v)$ at $(z_0;0,\kappa)$, for any $\kappa$ in $\mathbb{C}$, and
$A(z;0,\kappa)=\varepsilon K$.  The existence of a solution of the appropriate form again follows from Painlev\'e's lemma (lemma \ref{regularsequence}).  In particular, we see that $z-z_0\sim \varepsilon u^K$, so $u\sim c_0(z-z_0)^{1/K}$, where $c_0^K=\varepsilon$.

\vskip 3mm

\noindent{\bf Case 3:} $A=0$.\\
In this case the segments of the curve $\gamma$ on which $y$ is small can be deformed through a region on which $y$ is analytic such that $y$ is bounded away from $0$ on the resulting curve.  The proof is essentially identical to Shimomura (\cite{shimomura:03}, lemma 2.2 and ``case (iii)'' following remark 3.1), see also \cite{shimomura:07}.

\vskip 5mm

\noindent{\bf Proof of theorem \ref{main}, part 2}
The fact that these series converge follows from the existence of the regular systems of equations for $u$ and $v$ (i.e., equations (\ref{evenueqn}--\ref{evenveqn}) in the even case and equations (\ref{oddueqn}--\ref{oddveqn}) in the odd $N$ case).
\hfill $\Box$

\section{Singularity clustering along infinite length paths \label{sect-accumulate}}

So far we have only considered singularities that can be reached by analytic continuation along finite length paths.
The following proof is identical to that of the corresponding result in Smith \cite{smith:53} concerning equation
\eqref{smith} with $\mbox{deg}\,f>\mbox{deg}\,g$.

\vskip 5mm

\noindent{\bf Proof of theorem \ref{main}, part 4}

 Let $y$ be any solution of equation
\eqref{original} that can be analytically continued along a curve $\gamma\subset \Omega$ up to but not including the point $z_0$.  Let $D$ be an open disk in $\Omega$ centred at $z_0$.   Without loss of generality, we assume that $\gamma\subset D$.  Let $z_1$ be any point on $\gamma$.  There exists a straight line segment $l\in D$ starting from  $z_1$ such that  analytic continuation of $y$ along $l$ ends at a singularity $z_2$ of $y$ in $D$.  If this were not the case then $y$ would be analytic throughout $D$, contradicting the assumption that $y$ is singular at $z_0$.  Since the singularity at $z_2$ can be reached by analytic continuation along a finite length curve $l$, it follows from part 3 of theorem \ref{main} that it is of the algebraic type described.
\hfill$\Box$

\vskip 5mm

Next we summarize Smith's analysis \cite{smith:53} of equation \eqref{smith}.  Although this example is not of the form of equation \eqref{original}, it shows that the phenomenon described in part 4 of theorem \ref{main} can indeed occur for some simple equations.  

Smith showed that for any algebraic singularity of a solution of equation \eqref{smith} at a point $z_0\in\mathbb{C}$, the function
$$
\Phi(z):=y'(z)+y(z)^4
$$
is bounded in a neighbourhood of $z_0$.
Furthermore, he showed that the only possible singularities of any solution of equation \eqref{smith} are algebraic singularities or points at which such singularities accumulate along infinite length curves as described above.  He shows that such an accumulation of algebraic singularities actually occurs by constructing a path $\gamma$ ending at a finite point in the complex plane such that $\Phi$ tends to infinity along $\gamma$.  

Smith found the following parametric representation of a one-parameter family of solutions of equation \eqref{smith}:
\begin{eqnarray}
y(z) &=& 2^{-1/3}3^{1/6}x^{1/6}\left[
\frac{J_{-2/3}(x)}{J_{1/3}(x)}
\right]^{1/2},
\label{ypara}\\
z &=& 3^{-1/2}\int_{x_0}^x \left[
\frac{J_{1/3}(\xi)}{J_{-2/3}(\xi)}
\right]^{1/2}
\frac{{\rm d}\xi}{\xi^{1/2}},
\end{eqnarray}
and
$$
x=\frac{4i}{3}\Phi(z)^{3/2},
$$
where $x_0$ is an arbitrary constant and $J_\nu$ is the Bessel function of order $\nu$.
 Smith showed that as $x$ traces the path
 $$
 x=t-{\rm i}\cos 2\left(t-\frac{\pi}{6}-\frac{\pi}{4}\right),\qquad 2\pi\le t<\infty,
 $$
 $z$ traces a path $\gamma$ ending at a finite point $z_*$ in the finite complex plane.  Since $\Phi(z)\to\infty$ along this path, $z_*$ must be an accumulation point of algebraic singularities.

The kind of accumulation of singularities in Smith's example is far more complicated than the accumulation of poles in Painlev\'e's example \eqref{painleve}.  The accumulation point of poles in Painlev\'e's example can be reached by analytic continuation along finite length curves, so the standard methods used in the previous section show that no such accumulation can occur in many equations.  

Also, the more complicated accumulation of singularities along infinite length curves described by Smith's example (and allowed for in part 4 of theorem \ref{main}) does not arise as a issue in the standard proofs that solutions of equations such as the first Painlev\'e equation \eqref{p1} are meromorphic.  In these proofs one assumes that a solution $y$ is analytic at a point $z_0$ and, arguing by contradiction, let $a$ be the location of the nearest singularity that is not a pole.  Since $a$ can be reached by analytic continuation along a finite length curve, the standard arguments show that $y$ is either regular or has a pole at $a$, i.e., we obtain a contradiction.   It is only when the finite distance singularities are branch points that the  accumulation in the finite plane along infinite length curves becomes a possibility.

\section{Discussion}
In this paper we have shown that when equation \eqref{original} admits ``enough'' formal algebraic series solutions at movable singularities, then such singularities are the only ones that can be reached by analytic continuation along curves of finite length.  The existence of these formal series solutions is equivalent to one ($N$ even) or two ($N$ odd) resonance conditions being satisfied.  When these conditions are not satisfied, it is straightforward to show that there are formal series solutions of equation \eqref{canonical} the form
$$
y(z)=\sum_{n=0}^\infty c_n(\log [z-z_0])\,(z-z_0)^{\frac{n-2}{N-1}},
$$
where $c_0$ is a constant satisfying $c_0^{N-1}=1$ and each $c_n$ is a polynomial.  If one tries to extend the proof above to include series of this form, one can construct a function $W$, which now includes a factor involving $\log y$.  Similarly, the corresponding $u$-$v$ systems (\ref{evenueqn}--\ref{evenveqn}) and (\ref{oddueqn}--\ref{oddveqn}) have terms involving $\log u$.   In this setting, it is not sufficient to make estimates involving $|y|$ alone as one also needs to control the argument of $y$.

Whether algebraic branch points in solutions of equation \eqref{original} can accumulate in the manner described by part 4 for theorem \ref{main} remains an open question.

\vskip 5mm

\noindent{\large\bf Acknowledgments}
\\
The research reported in this paper was supported by a project grant from the Leverhulme Trust.  
The second author was supported through and EPSRC Advanced Research Fellowship.
We also acknowledge the support of the European Commission's Framework 6 ENIGMA Network and the
European Science Foundation's MISGAM Network.

\bibliographystyle{plain}

\end{document}